\documentclass[reqno]{amsart}

\newtheorem{Theorem}{Theorem}[section]

\theoremstyle{remark}
\newtheorem{Remark}[Theorem]{Remark}
\numberwithin{equation}{section}

\allowdisplaybreaks

\begin{document}

\title[Another proof of the ${}_6\psi_6$ sum]
{Another proof of Bailey's
$\hbox{${}_{\boldsymbol 6}\boldsymbol{\psi_{\boldsymbol 6}}$}$ summation}

\author[Fr\'ed\'eric Jouhet]{Fr\'ed\'eric Jouhet$^*$}
\address{Institut Girard Desargues, Universit\'e Claude Bernard (Lyon 1),
69622 Villeurbanne Cedex, France}
\email{jouhet@euler.univ-lyon1.fr}
\urladdr{http://igd.univ-lyon1.fr/home/jouhet}

\author[Michael Schlosser]{Michael Schlosser$^{**}$}
\address{Institut f\"ur Mathematik der Universit\"at Wien,
Strudlhofgasse 4, A-1090 Wien, Austria}
\email{schlosse@ap.univie.ac.at}
\urladdr{http://www.mat.univie.ac.at/{\textasciitilde}schlosse}

\thanks{$^*$Fully supported by EC's IHRP Programme,
grant HPRN-CT-2001-00272, ``Algebraic Combinatorics in Europe'', while
visiting the University of Vienna from February to August 2003.}
\thanks{$^{**}$Fully supported by an APART fellowship of
the Austrian Academy of Sciences}
\date{December 11, 2003}

\begin{abstract}
Adapting a method used by Cauchy, Bailey, Slater, and more recently,
the second author, we give a new proof of Bailey's celebrated
${}_6\psi_6$ summation formula.
\end{abstract}

\maketitle

\section{Introduction}

In \cite{Sc}, one of the authors presented a new proof of
Ramanujan's ${}_1\psi_1$ summation formula
(cf.\ \cite[Appendix~(II.29)]{GR}),
\begin{equation}\label{1psi1}
{}_1\psi_1\!\left[\begin{matrix}a\\b\end{matrix};q,z\right]=
\frac{(q,b/a,az,q/az)_\infty}{(b,q/a,z,b/az)_\infty}
\end{equation}
(the notation is defined at the end of this introduction),
valid for $|q|<1$ and $|b/a|<|z|<1$.
This proof used a standard method
for obtaining a bilateral identity from a unilateral terminating identity,
a method already utilized by Cauchy~\cite{Ca}
in his second proof of Jacobi's~\cite{Jac} famous triple product identity
(see \eqref{jtpi}), a special case of Ramanujan's formula \eqref{1psi1}.

The same method (which is referred to as ``Cauchy's method'' in the sequel)
had also been exploited by Bailey~\cite[Secs.~3 and 6]{Ba1},
\cite{Ba2}, and Slater~\cite[Sec.~6.2]{Sl}.

It was conjectured in \cite[Remark~3.2]{Sc} that {\em any} bilateral
sum can be obtained from an appropriately chosen terminating identity
by Cauchy's method (without appealing to analytic continuation).
However, at the same place it was also pointed out that
it was already not known whether Bailey's~\cite[Eq.~(4.7)]{Ba1}
very-well-poised $_6\psi_6$ summation
(cf.~\cite[Appendix~(II.33)]{GR}),
\begin{multline}\label{6psi6}
{}_6\psi_6\!\left[\begin{matrix}q\sqrt{a},-q\sqrt{a},b,c,d,e\\
\sqrt{a},-\sqrt{a},aq/b,aq/c,aq/d,aq/e\end{matrix};
q,\frac{qa^2}{bcde}\right]\\
=\frac{(q,aq,q/a,aq/bc,aq/bd,aq/be,aq/cd,aq/ce,aq/de)_\infty}
{(q/b,q/c,q/d,q/e,aq/b,aq/c,aq/d,aq/e,a^2q/bcde)_\infty}
\end{multline}
(again, see the end of this introduction for the notation),
where $|q|<1$ and $|a^2q/bcde|<1$,
would follow from such an identity. It is maybe interesting to mention
that Bailey's ${}_6\psi_6$ summation \eqref{6psi6},
although it contains more parameters
than Ramanujan's ${}_1\psi_1$ summation \eqref{1psi1},
does {\em not} include the latter as a special case.

While the conjecture of \cite[Remark~3.2]{Sc} remains open, this paper
features a new derivation of Bailey's ${}_6\psi_6$ summation formula
using a variant of Cauchy's method. 
After explaining some notation in the end of this introduction,
the proof from \cite[Sec.~3]{Sc}
of Ramanujan's ${}_1\psi_1$ summation
is being reviewed in Section~\ref{seccm}, for illustration.
The starting point for the derivation
of Bailey's ${}_6\psi_6$ summation in Section~\ref{secmain}
is Bailey's terminating very-well-poised
${}_{10}\phi_9$ transformation \eqref{Ba}. After appropriately applying
Cauchy's method to both sides of this identity, a specific transformation
of ${}_6\psi_6$ series is obtained. The resulting transformation
is then iterated, after which a suitable specialization yields the desired
${}_6\psi_6$ summation.

Other proofs of Bailey's very-well-poised ${}_6\psi_6$ summation
have been given by Bailey~\cite{Ba1}, Slater and Lakin~\cite{SL},
Andrews~\cite{And}, Askey and Ismail~\cite{AI},
Askey~\cite{AII}, and the second author~\cite{Sc6psi6}.

The authors plan to give an account of Cauchy's method
applied to terminating {\em quadratic}, {\em cubic} and {\em quartic}
identities~\cite[Sec.~3.8]{GR} elsewhere.

{\bf Notation:}
It is appropriate to recall some standard notation for \emph{$q$-series} and
\emph{basic hypergeometric series}.

Let $q$ be a fixed complex parameter (the ``base'') with $0<|q|<1$.
The \emph{$q$-shifted factorial} is defined for any complex
parameter $a$ by
\begin{equation*}
(a)_\infty\equiv (a;q)_\infty:=\prod_{j\geq 0}(1-aq^j)
\end{equation*}
and
\begin{equation*}
(a)_k\equiv (a;q)_k:=\frac{(a;q)_\infty}{(aq^k;q)_\infty},
\end{equation*}
where $k$ is any integer.
Since the same base $q$ is used throughout this paper,
it may be readily omitted (in notation) which will not lead
to any confusion. For brevity, write
\begin{equation*}
(a_1,\cdots,a_m)_k:=(a_1)_k\cdots(a_m)_k,
\end{equation*}
where $k$ is an integer or infinity.
Further, recall the definition of \emph{basic hypergeometric series},
\begin{equation*}
{}_s\phi_{s-1}\!\left[\begin{matrix}a_1,\dots,a_s\\
b_1,\dots,b_{s-1}\end{matrix};q,z\right]:=
\sum_{k=0}^\infty\frac{(a_1,\dots,a_s)_k}{(q,b_1,\dots,b_{s-1})_k}z^k,
\end{equation*}
and of \emph{bilateral basic hypergeometric series},
\begin{equation*}
{}_s\psi_s\!\left[\begin{matrix}a_1,\dots,a_s\\
b_1,\dots,b_s\end{matrix};q,z\right]:=
\sum_{k=-\infty}^\infty\frac{(a_1,\dots,a_s)_k}{(b_1,\dots,b_s)_k}z^k.
\end{equation*}

See Gasper and Rahman's text~\cite{GR} for a
comprehensive study on the theory of basic hypergeometric series.
In particular, the computations in this paper rely on some elementary
identities for $q$-shifted factorials, listed in \cite[Appendix~I]{GR}.

\section{Cauchy's method and Ramanujan's ${}_1\psi_1$ summation}\label{seccm}

For illustration of ``Cauchy's method'', it is convenient to
review the second author's proof~\cite[Sec.~3]{Sc} of
Ramanujan's ${}_1\psi_1$ summation. A closely related analysis is
being applied in Section~\ref{secmain} to prove Bailey's
very-well-poised ${}_6\psi_6$ summation \eqref{6psi6}.

In Jackson's~\cite{Ja32} $q$-Pfaff--Saalsch\"utz summation
(cf.\ \cite[Appendix~(II.12)]{GR}),
\begin{equation}\label{qps}
{}_3\phi_2\!\left[\begin{matrix}a,b,q^{-n}\\
c,abq^{1-n}/c\end{matrix};q,q\right]=\frac{(c/a,c/b)_n}{(c,c/ab)_n},
\end{equation}
first replace $n$ by $2n$ and then shift the index of summation
by $n$ such that the new sum runs from $-n$ to $n$:
\begin{multline*}
\frac{(c/a,c/b)_{2n}}{(c,c/ab)_{2n}}=
\sum_{k=0}^{2n}\frac{(a,b,q^{-2n})_k}{(q,c,abq^{1-2n}/c)_k}q^k\\
=\frac{(a,b,q^{-2n})_n}{(q,c,abq^{1-2n}/c)_n}q^n
\sum_{k=-n}^{n}\frac{(aq^n,bq^n,q^{-n})_k}{(q^{1+n},cq^n,abq^{1-n}/c)_k}q^k.
\end{multline*}
Next, replace $a$ by $aq^{-n}$, and $c$ by $cq^{-n}$, respectively, to get
\begin{multline*}
\sum_{k=-n}^{n}\frac{(a,bq^n,q^{-n})_k}{(q^{1+n},c,abq^{1-n}/c)_k}q^k=
\frac{(c/a,cq^{-n}/b)_{2n}}{(cq^{-n},c/ab)_{2n}}
\frac{(q,cq^{-n},abq^{1-2n}/c)_n}{(aq^{-n},b,q^{-2n})_n}q^{-n}\\
=\frac{(c/a)_{2n}}{(q)_{2n}}\frac{(q,q,c/b,bq/c)_n}{(c,q/a,b,c/ab)_n}.
\end{multline*}
Now, one may let $n\to\infty$, assuming $|c/ab|<1$ and $|b|<1$,
while appealing to Tannery's theorem~\cite{Br} for being allowed to
interchange limit and summation. This gives
\begin{equation*}
\sum_{k=-\infty}^\infty\frac{(a)_k}{(c)_k}\left(\frac{c}{ab}\right)^k=
\frac{(q,c/a,c/b,bq/c)_\infty}{(c,q/a,b,c/ab)_\infty},
\end{equation*}
where $|c/a|<|c/ab|<1$. Finally, replacing $b$ by $c/az$ and then $c$ by $b$
gives \eqref{1psi1}.

\begin{Remark}
Cauchy~\cite{Ca} had applied the above method to a special case
of \eqref{qps}. (After all, \eqref{qps} was not available to him yet.)
His starting point was the terminating $q$-binomial theorem,
\begin{equation}\label{qbin}
{}_1\phi_0\!\left[\begin{matrix}q^{-n}\\
-\end{matrix}\,;q,zq^n\right]=(z)_n,
\end{equation}
an identity already known to Euler~\cite{Eu},
which can be derived from \eqref{qps} by first doing the substitution
$b\mapsto c/z$, then $c\mapsto b$, and then letting
$a\to\infty$ and $b\to 0$. As a result of ``bilateralizing'' \eqref{qbin}
by the above procedure, Cauchy recovered Jacobi's~\cite{Jac}
well-known triple product identity,
\begin{equation}\label{jtpi}
\sum_{k=-\infty}^{\infty}(-1)^kq^{\binom k2}z^k=
(q,z,q/z)_{\infty}.
\end{equation}
Note that \eqref{jtpi} can be obtained from \eqref{1psi1} by first
replacing $z$ by $z/a$, and then letting $a\to\infty$ and $b\to 0$.
\end{Remark}

\section{Proof of Bailey's ${}_6\psi_6$ summation formula}\label{secmain}

In order to prove Bailey's ${}_6\psi_6$ summation formula \eqref{6psi6}
by Cauchy's method, one should start with a known terminating
identity which contains enough parameters. If one considers
Jackson's~\cite{Ja87} summation formula~\cite[Appendix~(II.22)]{GR},
\begin{multline}\label{Ja}
{}_{8}\phi_7\!\left[\begin{matrix}a,\,q\sqrt{a},-q\sqrt{a},
b,c,d,a^2q^{n+1}/bcd,q^{-n}\\
\sqrt{a},-\sqrt{a},aq/b,aq/c,aq/d,bcdq^{-n}/a,aq^{n+1}\end{matrix};
q,q\right]\\
=\frac{(aq,aq/bc,aq/bd,aq/cd)_n}{(aq/b,aq/c,aq/d,aq/bcd)_n},
\end{multline}
it becomes apparent that at least one parameter is missing here
for the given purpose. As a next step, one might consider
Watson's transformation formula of a terminating
${}_8\phi_7$ into a multiple of a ${}_4\phi_3$~\cite[Appendix~(III.18)]{GR},
which involves an additional parameter, and apply Cauchy's method
(to both sides of the transformation).
In fact, this was undertaken by Bailey~\cite{Ba2} who obtained
by this procedure, and some symmetry argument,
a transformation for ${}_2\psi_2$ series (see also \cite[Ex.~5.11]{GR}).

The next level object in the hierarchy of identities for (very-well-poised)
basic hypergeometric series is Bailey's~\cite{Ba0} transformation
formula~\cite[Appendix~(III.28)]{GR},
\begin{multline}\label{Ba}
{}_{10}\phi_9\!\left[\begin{matrix}a,\,q\sqrt{a},-q\sqrt{a},b,c,d,e,f,
\lambda aq^{n+1}/ef,q^{-n}\\
\sqrt{a},-\sqrt{a},aq/b,aq/c,aq/d,aq/e,aq/f,
efq^{-n}/\lambda,aq^{n+1}\end{matrix};q,q\right]\\
=\frac{(aq,aq/ef,\lambda q/e,\lambda q/f)_n}
{(aq/e,aq/f,\lambda q/ef,\lambda q)_n}\\\times
{}_{10}\phi_9\!\left[\begin{matrix}\lambda,\,
q\sqrt{\lambda},-q\sqrt{\lambda},\lambda
b/a,\lambda c/a,\lambda d/a,e,f,\lambda aq^{n+1}/ef,q^{-n}\\
\sqrt{\lambda},-\sqrt{\lambda},aq/b,aq/c,aq/d,\lambda q/e,
\lambda q/f,efq^{-n}/a,\lambda q^{n+1}\end{matrix};q,q\right],
\end{multline}
where $\lambda=qa^2/bcd$. A standard proof of \eqref{Ba}
involves two applications of \eqref{Ja}, together with an
interchange of summations, see \cite[Sec.~2.9]{GR}. Jackson's summation
\eqref{Ja} itself can be proved in various ways, see e.~g.\
Slater~\cite[Sec.~3.3.1]{Sl}, or Gasper and Rahman~\cite[Sec.~2.6]{GR}.

Note that for
$b$, $c$ or $d\to\infty$, \eqref{Ba} reduces to Watson's transformation
formula. By a further specialization one obtains Jackson's
${}_8\phi_7$ summation formula which may also be derived 
directly from \eqref{Ba} by letting $b=aq/c$ (thus $\lambda=a/d$).

To prove \eqref{6psi6}, start from \eqref{Ba}.
First, replace $n$ by $2n$, shift the index of summation by $n$ on
both sides, and obtain
\begin{multline*}
\frac{(1-aq^{2n})}{(1-a)}
\frac{(a,b,c,d,e,f,\lambda aq^{2n+1}/ef,q^{-2n})_n}
{(q,aq/b,aq/c,aq/d,aq/e,aq/f,efq^{-2n}/\lambda,aq^{2n+1})_n}q^n\\\times
\sum_{k=-n}^n\frac{(1-aq^{2n+2k})}{(1-aq^{2n})}
\frac{(aq^n,bq^n,cq^n,dq^n)_k}
{(q^{1+n},aq^{1+n}/b,aq^{1+n}/c,aq^{1+n}/d)_k}\\\times
\frac{(eq^n,fq^n,\lambda aq^{3n+1}/ef,q^{-n})_k}
{(aq^{1+n}/e,aq^{1+n}/f,efq^{-n}/\lambda,aq^{3n+1})_k}q^k\\
=\frac{(aq,aq/ef,\lambda q/e,\lambda q/f)_{2n}}
{(aq/e,aq/f,\lambda q/ef,\lambda q)_{2n}}\\\times
\frac{(1-\lambda q^{2n})}{(1-\lambda)}\frac{(\lambda,\lambda b/a,
\lambda c/a,\lambda d/a,e,f,\lambda aq^{2n+1}/ef,q^{-2n})_n}
{(q,aq/b,aq/c,aq/d,\lambda q/e,\lambda q/f,efq^{-2n}/a,
\lambda q^{2n+1})_n}q^n\\\times
\sum_{k=-n}^n\frac{(1-\lambda q^{2n+2k})}{(1-\lambda q^{2n})}
\frac{(\lambda q^n,\lambda bq^n/a,\lambda cq^n/a,\lambda dq^n/a)_k}
{(q^{1+n},aq^{1+n}/b,aq^{1+n}/c,aq^{1+n}/d)_k}\\\times
\frac{(eq^n,fq^n,\lambda aq^{3n+1}/ef,q^{-n})_k}
{(\lambda q^{1+n}/e,\lambda q^{1+n}/f,efq^{-n}/a,\lambda q^{3n+1})_k}q^k.
\end{multline*}
Next, replace $a$, $c$, $d$, $e$ and $f$ by
$aq^{-2n}$, $cq^{-n}$, $dq^{-n}$, $eq^{-n}$ and $fq^{-n}$, respectively.
This gives
\begin{multline*}
\sum_{k=-n}^n\frac{(1-aq^{2k})}{(1-a)}
\frac{(aq^{-n},bq^n,c,d,e,f,\lambda aq^{n+1}/ef,q^{-n})_k}
{(q^{1+n},aq^{1-n}/b,aq/c,aq/d,aq/e,aq/f,efq^{-n}/\lambda,aq^{n+1})_k}q^k\\
=\frac{(1-aq^{-2n})}{(1-a)}
\frac{(aq^{1-n}/e,aq^{1-n}/f,efq^{-2n}/\lambda,aq)_n}
{(aq^{-2n},b,cq^{-n},dq^{-n})_n}\\\times
\frac{(aq^{1-2n},aq/ef,\lambda q^{1-n}/e,\lambda q^{1-n}/f)_{2n}}
{(aq^{1-n}/e,aq^{1-n}/f,\lambda q/ef,\lambda q^{1-2n})_{2n}}\\\times
\frac{(1-\lambda)}{(1-\lambda q^{-2n})}
\frac{(\lambda q^{-2n},\lambda b/a,\lambda cq^{-n}/a,\lambda dq^{-n}/a)_n}
{(\lambda q^{1-n}/e,\lambda q^{1-n}/f,efq^{-2n}/a,\lambda q)_n}\\\times
\sum_{k=-n}^n\frac{(1-\lambda q^{2k})}{(1-\lambda)}
\frac{(\lambda q^{-n},\lambda bq^{n}/a,
\lambda c/a,\lambda d/a,e,f,\lambda aq^{n+1}/ef,q^{-n})_k}
{(q^{1+n},aq^{1-n}/b,aq/c,aq/d,
\lambda q/e,\lambda q/f,efq^{-n}/a,\lambda q^{n+1})_k}q^k\\
=\frac{(\lambda q/e,\lambda q/f,aq,\lambda b/a,aq/\lambda c,
aq/\lambda d,q/a,aq/ef)_n}{(aq/e,aq/f,b,\lambda q,q/c,q/d,q/\lambda,
\lambda q/ef)_n}\\\times
\sum_{k=-n}^n\frac{(1-\lambda q^{2k})}{(1-\lambda)}
\frac{(\lambda q^{-n},\lambda bq^{n}/a,
\lambda c/a,\lambda d/a,e,f,\lambda aq^{n+1}/ef,q^{-n})_k}
{(q^{1+n},aq^{1-n}/b,aq/c,aq/d,
\lambda q/e,\lambda q/f,efq^{-n}/a,\lambda q^{n+1})_k}q^k.
\end{multline*}
Now, one may let $n\to\infty$, assuming $|qa^2/cdef|<1$ while appealing to
Tannery's theorem, which yields the following transformation formula:
\begin{multline}\label{b2}
{}_6\psi_6\!\left[\begin{matrix}q\sqrt{a},-q\sqrt{a},c,d,e,f\\
\sqrt{a},-\sqrt{a},aq/c,aq/d,aq/e,aq/f\end{matrix};
q,\frac{qa^2}{cdef}\right]\\
=\frac{(aq,q/a,aq/ef,aq/cd,\lambda q/e,\lambda
 q/f,aq/\lambda c,aq/\lambda d)_\infty}{(aq/e,aq/f,q/c,q/d,\lambda
 q,q/\lambda,\lambda q/ef,b)_\infty}\\\times
{}_6\psi_6\!\left[\begin{matrix}q\sqrt{\lambda},-q\sqrt{\lambda},
\lambda c/a,\lambda d/a,e,f\\
\sqrt{\lambda},-\sqrt{\lambda},aq/c,aq/d,\lambda
q/e,\lambda q/f\end{matrix};q,\frac{qa^2}{cdef}\right],
\end{multline}
where $\lambda=qa^2/bcd$.
In this identity, the right-hand side involves one more parameter (namely $b$)
than the left-hand side. Note that when $b=aq/cd$ (whence $\lambda=a$),
this identity is trivial. It should also be pointed out that
\eqref{b2} is {\em not} Slater's~\cite{Sl0} transformation
\cite[Eq.~(5.5.3)]{GR}, the latter involving more symmetric parameters.

A possibility, of course, would be to directly specialize the extra
parameter $b$ in \eqref{b2} such that the ${}_6\psi_6$ on the
right-hand side reduces to a ${}_6\phi_5$ series, which can be summed
by the $n\to\infty$ case of \eqref{Ja}. However, this was not the idea
of the above derivaton of the ${}_6\psi_6$ transformation in \eqref{b2}.
Indeed, several proofs of Bailey's ${}_6\psi_6$ summation which 
make use of the nonterminating ${}_6\phi_5$ summation already exist,
see e.~g.\ Slater and Lakin~\cite{SL}, Andrews~\cite{And},
Askey and Ismail~\cite{AI}, and the second author~\cite{Sc6psi6}.

The clue is to {\em iterate} formula \eqref{b2}, more precisely, to
apply the same transformation to the ${}_6\psi_6$ on the right-hand
side of \eqref{b2} again, with the parameters $a$, $c$, $d$ and $e$,
respectively, being replaced by $\lambda$, $\lambda c/a$, $e$ and
$\lambda d/a$.
By this iteration an additional free parameter\footnote{On the contrary,
iteration of \cite[Eq.~(5.5.3)]{GR} would only yield an additional
{\em redundant} parameter and essentially reduce to the same identity.},
say $b'$, is obtained on the right-hand side.
The result is the following.
\begin{multline*}
{}_6\psi_6\!\left[\begin{matrix}q\sqrt{a},-q\sqrt{a},c,d,e,f\\
\sqrt{a},-\sqrt{a},aq/c,aq/d,aq/e,aq/f\end{matrix};
q,\frac{qa^2}{cdef}\right]\\
=\frac{(aq,q/a,aq/ef,aq/cd,\lambda q/e,aq/\lambda d)_\infty}
{(aq/e,aq/f,q/c,q/d,\lambda q/ef,b)_\infty}\\\times
\frac{(aq/df,aq/ec,a\lambda'q/\lambda
d,\lambda'q/f,aq/\lambda'c,\lambda q/\lambda'e)_\infty}{(aq/d,q/e,\lambda'
q,q/\lambda',a\lambda'q/\lambda df,b')_\infty}\\
\times {}_6\psi_6\!\left[\begin{matrix}q\sqrt{\lambda'},-q\sqrt{\lambda'},
\lambda' c/a,\lambda'e/\lambda,\lambda d/a,f\\
\sqrt{\lambda'},-\sqrt{\lambda'},aq/c,\lambda q/e,
\lambda'aq/\lambda d,\lambda' q/f\end{matrix};q,\frac{qa^2}{cdef}\right],
\end{multline*}
where  $\lambda=qa^2/bcd$ and $\lambda'=aq\lambda/b'ce$.

Now there are two extra parameters appearing on the right-hand side.
Take $\lambda=e$ (thus $b=qa^2/cde$), and $\lambda'=a/c$ (thus $b'=q$),
such that the last ${}_6\psi_6$, getting terminated from above and from below,
equals $1$, which immediately establishes Bailey's formula
\eqref{6psi6} (where $b$ has been replaced by $f$).

\begin{Remark}
In the above analysis, after applying Cauchy's method to both sides
of Bailey's ${}_{10}\phi_9$ transformation \eqref{Ba}, the resulting
${}_6\psi_6$ transformation \eqref{b2} first was iterated and then
specialized.
A natural question is what happens if one would start with the iterate
of \eqref{Ba}, listed in \cite[Ex.~2.19]{GR}, and then apply Cauchy's method.
In fact, the whole analysis would be similar to the one above.
In particular, one would also obtain a transformation of ${}_6\psi_6$
series (similar to but different from \eqref {b2} though) with an extra
free parameter on the right-hand side. Again, one can iterate the
transformation to obtain a second additional parameter, and then specialize
the two extra parameters such that the ${}_6\psi_6$ on the right-hand side
reduces to a sum of one single term only.
\end{Remark}

\end{document}